\documentclass[reqno,12pt,a4]{amsart}
\def\Sp{\mathop{\rm Sp}\nolimits}

\def\e{{\varepsilon}}

\def\g{{\gamma}}
\def\l{{\lambda}}

\def\D{{\Delta}}
\def\G{{\Gamma}}

\def\be{\begin{equation}}
\def\ee{\end{equation}}
\def\bea{\begin{eqnarray*}}
\def\eea{\end{eqnarray*}}

\newtheorem{thm}{Theorem}%[section]
\newtheorem{lem}[thm]{Lemma}
\newtheorem{prop}[thm]{Proposition}
\newtheorem{dfn}[thm]{Definition}

\newtheorem{cor}[thm]{Corollary}

\newcommand{\norm}[1]{\left\| #1 \right\|}

\def\Hp{{\mathcal {H}_{\pi}}}

\date{}
\title{On Almost Representations of Property (T) Groups}

\author{V. M. Manuilov}

\address{%
V. M. Manuilov\\Dept. of Mechanics and Mathematics\\ Moscow State University\\
Leninskie Gory, Moscow\\ 119992 Russia\\and\\Dept. of Mathematics\\
Harbin Institute of Technology\\
Harbin, 150001, P.R.C.}

\email{manuilov@mech.math.msu.su}

\thanks{The first named author was partially supported
by RFFI grant 05-01-00923.}
\author{Chao You}
\address{Chao You\\Dept. of Mathematics\\
Harbin Institute of Technology\\
Harbin, 150001, P.R.C.}

\email{hityou1982@hotmail.com}

\begin{document}
\maketitle

 \begin{abstract}
Property (T) for groups means a dichotomy: a representation either
has an invariant vector or all vectors are far from being
invariant. We show that, under a stronger condition of \.{Z}uk, a
similar dichotomy holds for almost representations as well.
 \end{abstract}

%%%%%%%%%%%%%%%%%%%%%%%%%%%%%%%%%%%%%%%%%%%%%%%%%%%%%%%%%%%%%%
\section{Introduction}
%%%%%%%%%%%%%%%%%%%%%%%%%%%%%%%%%%%%%%%%%%%%%%%%%%%%%%%%%%%%%%
Let $\G$ be a group generated by a symmetric finite set $S$ and
let $\pi:\G \rightarrow \mathcal {U}(\Hp)$ be a unitary
representation of $\G$. Suppose that $\G$ has the property (T) of
Kazhdan (i.e. the trivial representation is isolated in the dual
space of $\G$). We refer to \cite{HV} for basic information about
property (T). It is well known \cite{HRV} that the spectrum of
$\pi(x)=\frac{1}{|S|}\sum_{s\in\G}\pi(s)$ has a gap near 1:
 $$
\Sp(\pi(x))\subset[-1,1-c]\cup\{1\},
 $$
where $c$ is the Kazhdan constant for $\G$ (with respect to $S$).
In terms of the group $C^*$-algebra, this means that we can apply
a continuous function $f$ such that $f(1)=1$ and $f(t)=0$ for any
$t\in[-1,1-c]$ to $x=\frac{1}{|S|}\sum_{g\in S} g\in C^*(\G)$ to
obtain the canonical projection $p=f(x)\in C^*(\G)$ corresponding
to the trivial representation \cite{V}.

Our aim is to generalize the above property for the case of almost
representations of $\G$. Recall that, for $\varepsilon\geq 0$, an
$\varepsilon$-almost representation $\pi$ of $\G$ (with respect to
the given set $S$ of generators) is the map $\pi:S\to \mathcal
U(\Hp)$ satisfying
 $$
\|\pi(s_1s_2)-\pi(s_1)\pi(s_2)\|\leq\varepsilon
 $$
for any $s_1,s_2,s_1s_2\in S$ and $\pi(s^{-1})=\pi(s)^*$ for any
$s\in S$. This definition appeared in \cite{CGM} and then (in a
slightly different form) in \cite{MM}. If $\varepsilon=0$ (in the
case, when $\G$ is finitely presented and $S$ is sufficiently big)
then a 0-almost representation obviously generates a genuine
representation of $\G$. It is known that for some applications it
suffices for $\pi$ to be defined on $S$ only instead of the whole
$\G$. Any small perturbation of a genuine representation is an
almost representation, but there exist almost representations that
are far from any genuine representation \cite{voi}. One should
distinguish almost representations from other group `almost'
notions, e.g. quasi-representations, almost actions etc.
\cite{Stern,F-M}, which are completely different.

If we have an asymptotic representation (i.e. a continuous family
of $\varepsilon_t$-almost representations
$(\pi_t)_{t\in[0,\infty)}$ with
$\lim_{t\to\infty}\varepsilon_t=0$) then it follows from the
theory of $C^*$-algebra asymptotic homomorphisms that the spectrum
of $\pi_t(x)$ has a gap for $t$ sufficiently great: there is a
continuous function $\alpha=\alpha(t)>0$ such that
$\lim_{t\to\infty}\alpha(t)=0$ and
$\Sp(\pi_t(x))\subset[-1,1-c+\alpha(t)]\cup[1-\alpha(t),1]$.
Unfortunately, if we are interested in a single almost
representation, it may be impossible to include an almost
representation into an asymptotic one \cite{M}, and we don't know
how to check existence of a spectral gap because there is no nice
formula for the projection $p$.

Nevertheless, there is a condition, which is only slightly
stronger than the property (T) and which provides a gap in
$\Sp(\pi(x))$ for an almost representation $\pi$. The importance
of this condition was discovered by A. \.{Z}uk \cite{Zuk}. Let us
recall his construction.

It is supposed that the neutral element doesn't belong to $S$. A
finite graph $L(S)$ is assigned to the set $S$ of generators as
follows: the set of vertices of $L(S)$ is $S$ and the set $T$ of
edges of $L(S)$ is the set of all pairs $\{(s,s'):s,s',s^{-1}s'\in
S\}$. By including some additional elements in $S$, one can assume
that the graph $L(S)$ is connected \cite{Zuk}. For a vertex $s\in
L(S)$, let $deg(s)$ denote its {\it degree}, i.e. the number of
edges adjacent to $s$. Let $\Delta$ be a discrete Laplace operator
acting on functions $f$ defined on vertices of $L(S)$ by
 \be \label{1.4}
\Delta f(s)= f(s)-\frac{1}{deg(s)}\sum_{s'\sim s}f(s'),
 \ee
where $s'\sim s$ means that the vertex $s'$ is adjacent to the
vertex $s$. Operator $\D$ is a non-negative, self-adjoint operator
on the (finitedimensional) Hilbert space $l^2(L(S), deg)$ and zero
is a simple eigenvalue of $\D$. Let $\l_1(L(S))$ be the smallest
non-zero eigenvalue of $\D$. We say that a group $\G$ with the
generating set $S$ satisfies the \.{Z}uk's condition if
$\l_1(L(S))>\frac{1}{2}$. One of the main results of \cite{Zuk}
claims that the \.{Z}uk's condition implies property (T) with the
Kazhdan constant
$c=\frac{2}{\sqrt{3}}\left(2-\frac{1}{\lambda_1(L(S))}\right)$. We
appreciate \.{Z}uk's approach because it allows to work with
almost representations as well. The main result of this paper is
the following statement:

 \begin{thm}\label{1}
Let $\G$, $S$ satisfy the \.{Z}uk's condition and let $c$ be as
above. There is a continuous function
$\alpha=\alpha(\varepsilon)\geq 0$ such that $\alpha(0)=0$ and
 $$
\Sp\Bigl(\frac{1}{|S|}\sum_{s\in S}\pi(s)\Bigr)\subset
[-1,1-c/2+\alpha(\varepsilon)]\cup[1-\alpha(\varepsilon),1]
 $$
for any $\varepsilon$-almost representation $\pi$.
 \end{thm}

 \begin{cor}
For any $\varepsilon$-almost representation $\pi$ there exists an
$(\varepsilon+6|S|\alpha(\varepsilon))$-almost representation
$\pi'$ such that $\|\pi'(s)-\pi(s)\|\leq 3|S|\alpha(\varepsilon)$
for any $s\in S$ and $\pi'=\tau\oplus\sigma$, where $\tau$ is a
trivial representation and $\sigma$ is an
$(\varepsilon+6|S|\alpha(\varepsilon))$-almost representation
satisfying
 \be\label{s1}
\Sp\Bigl(\frac{1}{|S|}\sum_{s\in S}\sigma(s)\Bigr)\subset
[-1,1-c/2+(1+3|S|)\alpha(\varepsilon)].
 \ee
 \end{cor}

 \begin{proof}
Let $H\subset \Hp$ be the range of the spectral projection of
$\frac{1}{|S|}\sum_{s\in S}\pi(s)$ corresponding to the set
$[1-c+\alpha(\varepsilon),1]$. Then
$\|\pi(s)\xi-\xi\|\leq|S|\alpha(\varepsilon)\|\xi\|$ for any $s\in
S$ and for any $\xi\in H$ and if we write $\pi(s)$ as a matrix
$\left(\begin{smallmatrix}A&B\\C&D\end{smallmatrix}\right)$ with
respect to the decomposition $H\oplus H^\perp$ then
$\|B\|\leq|S|\alpha(\varepsilon)$ and
$\|C\|\leq|S|\alpha(\varepsilon)$, hence there exists a unitary
$D'$ such that $\|D'-D\|\leq 2|S|\alpha(\varepsilon)$. Put
$\pi'(s)=\left(\begin{smallmatrix}1&0\\0&D'\end{smallmatrix}\right)$.
Then $\|\pi'(s)-\pi(s)\|\leq 3|S|\alpha(\varepsilon)$ and $\pi'$
is obviously an $(\varepsilon+6|S|\alpha(\varepsilon))$-almost
representation, which is trivial on $H$. Hence $H^\perp$ is
$\pi(s)$-invariant for all $s\in S$. Since the restriction of
$\frac{1}{|S|}\sum_{s\in S}\pi(s)$ onto $H^\perp$ satisfies
$\bigl(\frac{1}{|S|}\sum_{s\in S}\pi(s)\bigr)|_{H^\perp}\leq
1-c/2+\alpha(\varepsilon)$, we get (\ref{s1}).
 \end{proof}

The remaining part of the paper is devoted to the proof of Theorem
\ref{1}. The proof follows the proof of \.{Z}uk for genuine
representations, but has additional argument because relations for
almost representations do not hold exactly, but only
approximately. It will be seen from the proof that one can take
$\alpha(\varepsilon)=O(\varepsilon^{2/5})$ in Theorem \ref{1}.

\section{Proof of the theorem}

The following Hilbert spaces and operators are defined exactly as
in \cite{Zuk}: It doesn't matter that $\pi$ is not a
representation here.
 \begin{dfn}[\cite{Zuk}]
{\rm For $r=0,1 \text{ and } 2$ let $C^r$ be the Hilbert spaces
defined as follows:
 \begin{eqnarray*}
C^0\!\!\!&=&\!\!\!\{u:u \in \Hp\};\  \langle u,w
\rangle_{C^0}=\langle u,w
\rangle_{\Hp}|T| \text{ for } u,w\in C^0;\\
C^1\!\!\!&=&\!\!\!\{f:S \rightarrow \Hp:f(s^{-1})=-\pi(s^{-1})f(s)
\text{ for all } s\in S \};\  \langle f,g \rangle_{C^1}=\sum_{s\in
S }\langle
f(s),g(s)\rangle_{\Hp}n(s);\\
C^2\!\!\!&=&\!\!\!\{g:T\rightarrow \Hp\};\  \langle f,g
\rangle_{C^2}=\sum_{t\in T }\langle f(t),g(t)\rangle_{\Hp},
 \end{eqnarray*}
where $n(s)=\#\{s'\in S:(s,s')\in T\}$. }
 \end{dfn}

Since the graph $L(S)$ is connected, $n(s)>0$ for every $s\in S$
and $n(s)=deg(s)$. Moreover, it is easy to see that
$n(s)=n(s^{-1})$ and $\sum_{s\in S}n(s)=|T|$.

\begin{dfn}[\cite{Zuk}]
{\rm Let us define linear operators $d_1:C^0\rightarrow C^1$ and
$d_2:C^1\rightarrow C^2$ as follows:
$$
d_1u(s)=\pi(s)u-u,\text{ for all } u\in C^0;
$$
$$
d_2f((s,s'))=f(s)-f(s')+\pi(s)f(s^{-1}s'), \text{ for all } f\in
C^1.
$$
}
\end{dfn}

 \begin{lem}[\cite{Zuk}]
One has $d^*_1f=-2\sum_{s\in S}f(s)\frac{n(s)}{|T|}$ for any $f\in
C^1$ and $\norm{d^*_1}_{C^1\rightarrow C^0}\leq 2$.
 \end{lem}

In \cite{Zuk} it is shown that $d_2d_1=0$ for any unitary
representation. However, if $\pi$ is only an almost representation
then this doesn't hold any more. One can only show that this
composition is small.

 \begin{lem}\label{5}
For any $u\in C^0$ and $(s,s')\in T$ one has \be \label{2.1}
\norm{d_2d_1u((s,s'))}_{\Hp}\leq \e \norm{u}_{\Hp} \ee
 \end{lem}
 \begin{proof} By the definitions of $d_1$ and $d_2$, we have
\begin{eqnarray*}
\norm{d_2d_1u((s,s'))}_\Hp&=&\norm{d_1u(s)-d_1u(s')+\pi(s)d_1u(s^{-1}s')}_\Hp\\
&=&\norm{(\pi(s)u-u)-(\pi(s')u-u)+\pi(s)(\pi(s^{-1}s')u-u)}_\Hp\\
&=&\norm{\pi(s')u-\pi(s)\pi(s^{-1}s')u}_\Hp \leq\e\norm{u}_\Hp,
\end{eqnarray*}
hence we have $\norm{d_2d_1u((s,s'))}_{\Hp}\leq \e \norm{u}_{\Hp}$.
 \end{proof}

 \begin{cor}
$\|d_2d_1\|_{C^2\to C^0}\leq\e$.
 \end{cor}

That's why we have to introduce two more (sub)spaces.

\begin{dfn}
{\rm For any $\beta\geq 0$ set
 $$
B^0(\beta)=\{P_\Omega(d^*_1d_1)(u):u\in C^0 \}\subset C^0,\qquad
B^1(\beta)=\overline{\{d_1u:u\in B^0(\beta)\}}\subset C^1,
 $$
where $P_\Omega$ is the spectral projection corresponding to
$\Omega=[\beta,+\infty)$.}
\end{dfn}

It is clear that $B^0(\beta)$ and $B^1(\beta)$ are invariant
subspaces for $d^*_1d_1$ and $d_1d^*_1$ respectively.

 \begin{prop} \label{tm1}
If there exists $c>0$ and $0<\delta<c/2$ such that for every $f\in
B^1(\frac{\delta^2}{|T|})$
 \be \label{2.4}
\langle d_1d^*_1f,f \rangle_{C^1}
>c\langle f, f \rangle_{C^1}
 \ee
then, for any $\e$-almost representation $\pi$, either there
exists $u\in C^0$ such that
 \be\label{6}
\norm{\pi(s)u-u}_\Hp<\delta\norm{u}_\Hp {\mbox{ \ for\ any\ }}
s\in S
 \ee
or
 \be \label{3}
\max_{s\in S}\norm{\pi(s)u-u}_\Hp\geq c/2\norm{u}_\Hp
 \ee
for every $u\in C^0$.
 \end{prop}
 \begin{proof}
First, we show that if there is no $u\in C^0$ satisfying (\ref{6})
then $B^0(\frac{\delta^2}{|T|})=C^0$. If this is not true then
there exists a non-zero vector $u^\bot$ orthogonal to
$B^0(\frac{\delta^2}{|T|})$. Since $\norm
{d^*_1d_1u^\bot}_{C^0}<\frac{\delta^2}{|T|}\norm{u^\bot}_{C^0}$,
we have $\langle d_1u^\bot,d_1u^\bot \rangle_{C^1}=\langle
u^\bot,d^*_1d_1u^\bot\rangle_{C^0}\leq\norm
{u^\bot}_{C^0}\norm{d^*_1d_1u^\bot}_{C^0}<
\frac{\delta^2}{|T|}\norm{u^\bot}^2_{C^0}$, which implies that
$\norm{d_1u^\bot}_{C^1}<\frac{\delta}{\sqrt{|T|}}\norm
{u^\bot}_{C^0}$. By definition of $\norm{\cdot}_{C^1}$, it is easy
to see that $\norm{\pi(s)u^\bot-u^\bot}_\Hp\leq
\norm{d_1u^\bot}_{C^1}<\frac{\delta}{\sqrt{|T|}}\sqrt{|T|}\norm
{u^\bot}_\Hp=\delta\norm{u^\bot}_\Hp$ for any $s\in S$, which is
in contradiction with the assumption.

Next we prove that (\ref{2.4}) implies that the operator
$d_1d^*_1:B^1(\frac{\delta^2}{|T|})\rightarrow
B^1(\frac{\delta^2}{|T|})$ has a bounded inverse. By (\ref{2.4}),
$d_1d^*_1(B^1(\frac{\delta^2}{|T|}))$ is closed in
$B^1(\frac{\delta^2}{|T|})$. If
$d_1d^*_1(B^1(\frac{\delta^2}{|T|}))$ were different from
$B^1(\frac{\delta^2}{|T|})$, there would exist a non-zero vector
$f \in B^1(\frac{\delta^2}{|T|})$ orthogonal to
$d_1d^*_1(B^1(\frac{\delta^2}{|T|}))$. Then we would have, by
(\ref{2.4}),
$$0=\langle f, d_1d^*_1(f)\rangle_{C^1}> c\langle f, f\rangle_{C^1}$$ which is a contradiction.

Thus $d_1d^*_1:B^1(\frac{\delta^2}{|T|})\rightarrow
B^1(\frac{\delta^2}{|T|})$ has a bounded inverse
$(d_1d^*_1)^{-1}:B^1(\frac{\delta^2}{|T|})\rightarrow
B^1(\frac{\delta^2}{|T|})$ and
$\norm{(d_1d^*_1)^{-1}}_{B^1(\frac{\delta^2}{|T|})\rightarrow
B^1(\frac{\delta^2}{|T|})}\leq c^{-1}$.

Now suppose that neither (\ref{6}) nor (\ref{3}) holds. Then there
is some $\gamma<c/2$ and some $u\in \Hp$ such that
$\norm{u}_\Hp=1$ and $\norm{\pi(s)u-u}_\Hp\leq \gamma$ for every
$s\in S$. Therefore
$$
\norm{d_1u}^2_{B^1(\frac{\delta^2}{|T|})}=\sum_{s\in
S}\norm{d_1u(s)}^2_\Hp n(s)=\sum_{s\in S}\norm{\pi(s)u-u}^2_\Hp
n(s)\leq \sum_{s\in S}\gamma^2n(s)=\gamma^2|T|
$$
which gives $\norm{d_1u}_{B^1(\frac{\delta^2}{|T|})}\leq \gamma
\sqrt{|T|}$. Then
\begin{eqnarray*}
\norm{d^*_1(d_1d^*_1)^{-1}d_1u}_{C^0}&\leq&\norm{d^*_1}_{B^1(\frac{\delta^2}{|T|})\rightarrow
C^0}
\cdot\norm{(d_1d^*_1)^{-1}}_{B^1(\frac{\delta^2}{|T|})\rightarrow B^1(\frac{\delta^2}{|T|})}\cdot\norm{d_1u}_{B^1(\frac{\delta^2}{|T|})}\\
&\leq&2\cdot c^{-1}\cdot \g\sqrt{|T|}<\sqrt{|T|}.
\end{eqnarray*}
By definition of the norm in $C^0$ one has then
$d^*_1(d_1d^*_1)^{-1}d_1u=u'$, whence $\norm{u'}_\Hp<1$, so the
vector $u-u'$ is non-zero. Finally,
$$
d_1(u-u')=d_1u-d_1(d^*_1(d_1d^*_1)^{-1})d_1u=d_1u-d_1u,
$$
which means that, for every $s\in S$,
$$
\pi(s)(u-u')-(u-u')=0.
$$
Thus $u-u'$ is a non-zero invariant vector and (\ref{6}) holds,
which gives a contradiction.
 \end{proof}

%%%%%%%%%%%%%%%%%%%%%%%%%%%%%%%%%%%%%%%%%%%%%%%%%%%%%%%%%%%%%%
%\subsection{Relation between the operators $d_2$ and $D$}
%%%%%%%%%%%%%%%%%%%%%%%%%%%%%%%%%%%%%%%%%%%%%%%%%%%%%%%%%%%%%%

Following \cite{Zuk}, define the operator $D:C^1\rightarrow C^2$
by
$$
Df((s_1, s_2))=f(s_1)-f(s_2),
$$
where $f\in C^1$ and $(s_1, s_2)\in T$.

In \cite{Zuk}, the relation between $d_2$ and $D$ was investigated
and it was shown that $\frac{1}{3}\langle d_2f,d_2f
\rangle_{C^2}=\langle D f,Df \rangle_{C^2}-\langle f,f
\rangle_{C^1}$. Since here an almost representation is engaged, we
have to estimate the difference between the left and the right
hand side.

 \begin{lem} For every $f\in C^1$ one has
 \be \label{2.5}
\langle f,f\rangle_{C^1}=\sum_{(s,s')\in T}\langle
f(s^{-1}s'),f(s^{-1}s')\rangle_\Hp;
 \ee
 \be \label{2.6}
\norm{d_2f((s,s'))+d_2f((s',s))}_\Hp\leq \e
\norm{f((s')^{-1}s)}_\Hp;
 \ee
 \be \label{2.7}
d_2f((s,s'))=-\pi(s)d_2f((s^{-1},s^{-1}s'));
 \ee
 \be \label{2.8}
\sum_{(s,s')\in T}\langle
d_2f((s^{-1},s^{-1}s')),-f(s^{-1}s')\rangle_\Hp=\sum_{(s,s')\in
T}\langle d_2f((s,s')),-f(s')\rangle_\Hp;
 \ee
 \be \label{2.9}
\norm{d_2f((s,s'))-\pi(s')d_2f(((s')^{-1},(s')^{-1}s))}_\Hp\leq\e\norm{f((s')^{-1}s)}_\Hp;
 \ee
 \be \label{2.10}
\Bigl|\sum_{(s,s')\in T}\langle d_2f((s, s')),
\pi(s)f(s^{-1}s')\rangle_\Hp-\frac{1}{3}\langle d_2f,
d_2f\rangle_{C^2}\Bigr|<\frac{5}{3}\e\norm{f}^2_{C^1}
 \ee
 \end{lem}
 \begin{proof}
The proof of (\ref{2.5}) and (\ref{2.7}) in \cite{Zuk} doesn't
depend on the property of $\pi$ being a representation.
\begin{eqnarray*}
d_2f((s,s'))&=&f(s)-f(s')+\pi(s)f(s^{-1}s')\\
&=&-(f(s')-f(s)+\pi(s)\pi(s^{-1}s')f((s')^{-1}s))\\
&=&-(f(s')-f(s)+\pi(s')f((s')^{-1}s)+\pi(s')f((s')^{-1}s)-\pi(s)\pi(s^{-1}s')f((s')^{-1}s))\\
&=&-d_2f((s',s))+(\pi(s')f((s')^{-1}s)-\pi(s)\pi(s^{-1}s')f((s')^{-1}s)),
\end{eqnarray*}
hence
\begin{eqnarray*}
\norm{d_2f((s,s'))+d_2f((s',s))}_\Hp&=&\norm{\pi(s')f((s')^{-1}s)-\pi(s)\pi(s^{-1}s')f((s')^{-1}s)}_\Hp\\
&\leq&\e \norm{f((s')^{-1}s)}_\Hp,
\end{eqnarray*}
which proves (\ref{2.6}).

By (\ref{2.6}) and (\ref{2.7}),
\begin{eqnarray*}
&&\!\!\!\!\!\!\!\!\!\!\!\!\!\!\!\!\!\!\!\!\!\norm{d_2f((s,s'))-\pi(s')d_2f(((s')^{-1},(s')^{-1}s))}_\Hp\\
&=&\|d_2f((s,s'))+d_2f((s',s))-d_2f((s',s))-\pi(s')d_2f(((s')^{-1},(s')^{-1}s))\|_\Hp\\
&\leq&\norm{d_2f((s,s'))+d_2f((s',s))}_\Hp \ \leq\
\e\norm{f((s')^{-1}s)}_\Hp,
\end{eqnarray*}
which proves (\ref{2.9}).

Consider the mapping $M:T\rightarrow T$,
$M((s,s'))=(s^{-1},s^{-1}s')=(t,t')$. Then it is easy to see $M$
is a bijection. Hence,
$$\sum_{(s,s')\in T}\langle
df((s^{-1},s^{-1}s')),-f(s^{-1}s')\rangle_\Hp=\sum_{(t,t')\in
T}\langle df((t,t')),-f(t')\rangle_\Hp,$$ which is just a matter
of notation. This proves (\ref{2.8}).

\begin{eqnarray*}
&&\sum_{(s,s')\in T}\langle d_2f((s, s')),
\pi(s)f(s^{-1}s')\rangle_\Hp\\&=&\frac{1}{3}\sum_{(s,s')\in
T}(\langle
d_2f((s, s')),\pi(s)f(s^{-1}s')\rangle_\Hp+\langle -\pi(s)d_2f((s^{-1},s^{-1}s')),\pi(s)f(s^{-1}s')\rangle_\Hp\allowdisplaybreaks\\
&&+\langle \pi(s')d_2f(((s')^{-1},(s')^{-1}s)),\pi(s)f(s^{-1}s')
\rangle_\Hp)+D_1\allowdisplaybreaks\\
&=&\frac{1}{3}\sum_{(s,s')\in T}(\langle d_2f((s, s')),
\pi(s)f(s^{-1}s')\rangle_\Hp+\langle
d_2f((s^{-1},s^{-1}s')),-f(s^{-1}s')\rangle_\Hp\allowdisplaybreaks\\&&+\langle
d_2f(((s')^{-1},(s')^{-1}s)),\pi((s')^{-1}s)f(s^{-1}s')
\rangle_\Hp)+D_1+D_2\allowdisplaybreaks\\
&=&\frac{1}{3}\sum_{(s,s')\in T}(\langle d_2f((s, s')),
\pi(s)f(s^{-1}s')\rangle_\Hp+\langle d_2f((s, s')), -f(s')\rangle_\Hp\allowdisplaybreaks\\
&&+\langle d_2f(((s')^{-1},(s')^{-1}s)),
-f((s')^{-1}s)\rangle_\Hp)+D_1+D_2\\
&=&\frac{1}{3}\sum_{(s,s')\in T}(\langle d_2f((s, s')),
\pi(s)f(s^{-1}s')\rangle_\Hp+\langle d_2f((s, s')), -f(s')\\
&&+\langle d_2f((s',s)),
-f(s)\rangle_\Hp)+D_1+D_2\\
&=&\frac{1}{3}\sum_{(s,s')\in T}(\langle d_2f((s, s')),
\pi(s)f(s^{-1}s')\rangle_\Hp+\langle d_2f((s, s')),
-f(s')\rangle_\Hp\\
&&+\langle d_2f((s, s')), f(s)\rangle_\Hp)
+D_1+D_2+D_3\\
&=&\frac{1}{3}\sum_{(s,s')\in T}(\langle d_2f((s, s')),f(s)-f(s')+
\pi(s)f(s^{-1}s')\rangle_\Hp+D_1+D_2+D_3\\
&=&\frac{1}{3}\sum_{(s,s')\in T}(\langle d_2f((s, s')),d_2f((s,
s'))\rangle_\Hp+D_1+D_2+D_3\\
&=&\frac{1}{3}\sum_{(s,s')\in T}\langle
d_2f,d_2f\rangle_{C^2}+D_1+D_2+D_3,
\end{eqnarray*}

where
\begin{eqnarray*}
&&D_1=\frac{1}{3}\sum_{(s,s')\in T}\langle d_2f((s,
s'))-\pi(s')d_2f(((s')^{-1},(s')^{-1}s)),\pi(s)f(s^{-1}s')\rangle_\Hp,\\
&&D_2=\frac{1}{3}\sum_{(s,s')\in T}\langle
d_2f(((s')^{-1},(s')^{-1}s)),(\pi((s')^{-1})\pi(s)-\pi((s')^{-1}s))f(s^{-1}s')\rangle_\Hp,\\
&&D_3=\frac{1}{3}\sum_{(s,s')\in T}\langle
d_2f((s,s'))+d_2f((s',s)),-f(s)\rangle_\Hp,
\end{eqnarray*}
hence $$\sum_{(s,s')\in T}\langle d_2f((s, s')),
\pi(s)f(s^{-1}s')\rangle_\Hp-\frac{1}{3}\sum_{(s,s')\in T}\langle
d_2f,d_2f\rangle_{C^2}=D_1+D_2+D_3.$$

By Cauchy inequality, definition of $\norm{\cdot}_{C^1}$ and
(\ref{2.5})---(\ref{2.9}), we have
\begin{eqnarray*}
|D_1|&=&\frac{1}{3}\Bigl|\sum_{(s,s')\in T}\langle d_2f((s,
s'))-\pi(s')d_2f(((s')^{-1},(s')^{-1}s)),\pi(s)f(s^{-1}s')\rangle_\Hp\Bigr|\allowdisplaybreaks\\
&\leq&\frac{1}{3}\sum_{(s,s')\in T}|\langle d_2f((s,
s'))-\pi(s')d_2f(((s')^{-1},(s')^{-1}s)),\pi(s)f(s^{-1}s')\rangle_\Hp|\allowdisplaybreaks\\
&\leq&\frac{1}{3}\sum_{(s,s')\in T}\norm{d_2f((s,
s'))-\pi(s')d_2f(((s')^{-1},(s')^{-1}s))}_\Hp\norm{\pi(s)f(s^{-1}s')}_\Hp\allowdisplaybreaks\\
&\leq&\frac{1}{3}\sum_{(s,s')\in
T}\e\norm{f((s')^{-1}s)}_\Hp\norm{f(s^{-1}s')}_\Hp\allowdisplaybreaks\\
&\leq&\frac{1}{3}\e\bigl(\sum_{(s,s')\in
T}\norm{f((s')^{-1}s)}^2_\Hp\bigr)^{1/2}\bigl(\sum_{(s,s')\in
T}\norm{f(s^{-1}s')}^2_\Hp\bigr)^{1/2}\ =\
\frac{1}{3}\e\norm{f}^2_{C^1},
\end{eqnarray*}

\begin{eqnarray*}
|D_2|&=&\frac{1}{3}\Bigl|\sum_{(s,s')\in T}\langle
d_2f(((s')^{-1},(s')^{-1}s)),(\pi((s')^{-1})\pi(s)-\pi((s')^{-1}s))f(s^{-1}s')\rangle_\Hp\Bigr|\\
&\leq&\frac{1}{3}\sum_{(s,s')\in T}\norm{
d_2f(((s')^{-1},(s')^{-1}s))}_\Hp\norm{(\pi((s')^{-1})\pi(s)-\pi((s')^{-1}s))f(s^{-1}s')}_\Hp\\
&\leq&\frac{1}{3}\sum_{(s,s')\in
T}\norm{f((s')^{-1})-f((s')^{-1}s)+\pi((s')^{-1})f(s)}_\Hp\e\norm{f(s^{-1}s')}_\Hp\\
&\leq&\frac{1}{3}\e\sum_{(s,s')\in
T}\norm{f((s')^{-1})}_\Hp\norm{f(s^{-1}s')}_\Hp+\frac{1}{3}\e\sum_{(s,s')\in
T}\norm{f((s')^{-1}s)}_\Hp\norm{f(s^{-1}s')}_\Hp\\
&&+\frac{1}{3}\e\sum_{(s,s')\in
T}\norm{f(s)}_\Hp\norm{f(s^{-1}s')}_\Hp\\
&\leq&\frac{1}{3}\e\bigl(\sum_{(s,s')\in
T}\norm{f((s')^{-1})}^2_\Hp\bigr)^{1/2}\bigl(\sum_{(s,s')\in
T}\norm{f(s^{-1}s')}^2_\Hp\bigr)^{1/2}\\
&&+\frac{1}{3}\e\bigl(\sum_{(s,s')\in
T}\norm{f((s')^{-1}s)}^2_\Hp\bigr)^{1/2}\bigl(\sum_{(s,s')\in
T}\norm{f(s^{-1}s')}^2_\Hp\bigr)^{1/2}\\
&&+\frac{1}{3}\e\bigl(\sum_{(s,s')\in
T}\norm{f(s)}^2_\Hp\bigr)^{1/2}\bigl(\sum_{(s,s')\in
T}\norm{f(s^{-1}s')}^2_\Hp\bigr)^{1/2}\\
&=&\frac{1}{3}\e\norm{f}^2_{C^1}+\frac{1}{3}\e\norm{f}^2_{C^1}+\frac{1}{3}\e\norm{f}^2_{C^1}
\ =\ \e\norm{f}^2_{C^1},
\end{eqnarray*}

\begin{eqnarray*}
|D_3|&=&\frac{1}{3}\Bigl|\sum_{(s,s')\in T}\langle
d_2f((s,s'))+d_2f((s',s)),-f(s)\rangle_\Hp\Bigr|\\
&\leq&\frac{1}{3}\sum_{(s,s')\in
T}\norm{d_2f((s,s'))+d_2f((s',s))}_\Hp\norm{f(s)}_\Hp\allowdisplaybreaks\\
&\leq&\frac{1}{3}\sum_{(s,s')\in
T}\e\norm{f((s')^{-1}s)}_\Hp\norm{f(s)}_\Hp\allowdisplaybreaks\\
&\leq&\frac{1}{3}\e\bigl(\sum_{(s,s')\in
T}\norm{f((s')^{-1}s)}^2_\Hp\bigr)^{1/2}\bigl(\sum_{(s,s')\in
T}\norm{f(s)}^2_\Hp\bigr)^{1/2} \ =\
\frac{1}{3}\e\norm{f}^2_{C^1}.\allowdisplaybreaks
\end{eqnarray*}

So
\begin{eqnarray*}
\Bigl|\sum_{(s,s')\in T}\langle d_2f((s, s')),
\pi(s)f(s^{-1}s')\rangle_\Hp-\frac{1}{3}\langle
d_2f,d_2f\rangle_{C^2}\Bigr|&=&|D_1+D_2+D_3|\\
&\leq&|D_1|+|D_2|+|D_3|\ \leq\ \frac{5}{3}\e\norm{f}^2_{C^1},
\end{eqnarray*}
which proves (\ref{2.10}).
 \end{proof}

 \begin{prop}\label{2.9}
For every $f\in C^1$ one has
$$
\bigl|\langle Df,Df \rangle_{C^2}-\frac{1}{3}\langle
d_2f,d_2f\rangle_{C^2}-\langle
f,f\rangle_{C^1}\bigr|\leq\frac{10}{3}\e\langle f,f\rangle_{C^1}
$$
 \end{prop}
 \begin{proof}
By definition of operator $D$, we have
\begin{eqnarray*}
\langle Df,Df \rangle_{C^2}&=&\sum_{(s,s')\in T}\langle
d_2f((s,s'))-\pi(s)f(s^{-1}s'),d_2f((s,s'))-\pi(s)f(s^{-1}s')\rangle_\Hp\\
&=&\sum_{(s,s')\in T}\langle
d_2f((s,s')),d_2f((s,s'))\rangle_\Hp+\sum_{(s,s')\in T}\langle
f(s^{-1}s'),f(s^{-1}s')\rangle_\Hp\\
&&-2\sum_{(s,s')\in T}\langle d_2f((s, s')),
\pi(s)f(s^{-1}s')\rangle_\Hp\\
&=&\langle d_2f,d_2f\rangle_{C^2}+\langle
f,f\rangle_{C^1}-\frac{2}{3}\langle d_2f,d_2f\rangle_{C^2}\\
&&+\frac{2}{3}\langle d_2f,d_2f\rangle_{C^2}-2\sum_{(s,s')\in
T}\langle d_2f((s, s')), \pi(s)f(s^{-1}s')\rangle_\Hp\\
&=&\frac{1}{3}\langle d_2f,d_2f\rangle_{C^2}+\langle
f,f\rangle_{C^1}\\
&&+2\bigl(\frac{1}{3}\langle
d_2f,d_2f\rangle_{C^2}-\sum_{(s,s')\in T}\langle d_2f((s, s')),
\pi(s)f(s^{-1}s')\rangle_\Hp\bigr),
\end{eqnarray*}
hence
\begin{eqnarray*}
&&\!\!\!\!\!\!\!\!\!\!\bigl|\langle Df,Df
\rangle_{C^2}-\frac{1}{3}\langle d_2f,d_2f\rangle_{C^2}-\langle
f,f\rangle_{C^1}\bigr|\\&=&2\Bigl|\frac{1}{3}\langle
d_2f,d_2f\rangle_{C^2}-\sum_{(s,s')\in T}\langle d_2f((s, s')),
\pi(s)f(s^{-1}s')\rangle_\Hp\Bigr|\ \leq\ \frac{10}{3}\e\langle
f,f\rangle_{C^1},
\end{eqnarray*}
which ends the proof of Proposition \ref{2.9}.
 \end{proof}

%%%%%%%%%%%%%%%%%%%%%%%%%%%%%%%%%%%%%%%%%%%%%%%%%%%%%%%%%%%%%%
%\subsection{Reduction to the graph $L(S)$}
%%%%%%%%%%%%%%%%%%%%%%%%%%%%%%%%%%%%%%%%%%%%%%%%%%%%%%%%%%%%%%
Note that every $f\in C^1$ can be considered as a function on
$L(S)$. It was shown in \cite{Zuk} (and the proof doesn't depend
on the property of $\pi$ to be a representation) that $\langle
f,f\rangle_{C^1}=\langle f,f\rangle_{L(S)}$,
 $$
\langle Df,Df\rangle_{C^2}=\sum_{(s,s')\in T}\langle
f(s)-f(s'),f(s)-f(s')\rangle_\Hp=2\langle \Delta f,f\rangle_{L(S)}
 $$
and
 \be\label{4}
\langle \Delta f,f\rangle_{L(S)}\geq\lambda_1(L(S))\langle
f,f\rangle_{C^1}-\frac{\lambda_1(L(S))}{4}\langle
d^*_1f,d^*_1f\rangle_{C^0}.
 \ee

From now on, for shortness' sake, we denote $\lambda_1(L(S))$ by
$\lambda_1$.
 \begin{lem}
For every $f\in C^1$ one has
 \be \label{L(S)} \frac{1}{3}\langle
d_2f,d_2f\rangle_{C^2}+\frac{1}{2}\lambda_1\langle
d^*_1f,d^*_1f\rangle_{C^0}\geq\bigl(2\lambda_1-1-\frac{10}{3}\e\bigr)\langle
f,f\rangle_{C^1}.
 \ee
 \end{lem}
 \begin{proof}
By (\ref{4}), we have
$$
2\lambda_1\langle f,f\rangle_{C^1}-\langle
f,f\rangle_{C^1}\leq2\langle \Delta f,f\rangle_{L(S)}-\langle
f,f\rangle_{C^1}+\frac{\lambda_1}{2}\langle
d^*_1f,d^*_1f\rangle_{C^0}.
$$
Since $\langle Df,Df\rangle_{C^2}=2\langle \Delta
f,f\rangle_{L(S)}$,
\begin{eqnarray*}
2\lambda_1\langle f,f\rangle_{C^1}-\langle f,f\rangle_{C^1}
&\leq&\langle Df,Df\rangle_{C^2}-\langle
f,f\rangle_{C^1}+\frac{\lambda_1}{2}\langle
d^*_1f,d^*_1f\rangle_{C^0}\\
&\leq&\frac{1}{3}\langle
d_2f,d_2f\rangle_{C^2}+\frac{10}{3}\e\langle
f,f\rangle_{C^1}+\frac{\lambda_1}{2}\langle
d^*_1f,d^*_1f\rangle_{C^0},
\end{eqnarray*}
which proves (\ref{L(S)}).
 \end{proof}

\begin{lem}
For any $f\in B^1(\frac{\delta^2}{|T|})$ one has
 $$
\norm{d_2f}^2_{C^2}\leq\frac{4|T|^2\e^2}{\delta^4}\norm{f}_{C^1}.
 $$
 \end{lem}
 \begin{proof}
By Lemma \ref{5}, for any $u\in B^0(\frac{\delta^2}{|T|})$ we have
 $$
\norm{d_2d_1u}^2_{C^2}=\sum_{(s,s')\in T}\norm{d_2d_1u((s,
s'))}^2_\Hp\leq\e^2\sum_{(s,s')\in T}\norm{u}^2_\Hp
=\e^2|T|\norm{u}^2_\Hp\\=\e^2\norm{u}^2_{C^0}.
 $$

Since $u\in B^0(\frac{\delta^2}{|T|})$, it follows that
$\frac{\delta^2}{|T|}\norm{u}_{C^0}\leq\norm{d^*_1d_1u}_{C^0}\leq2\norm{d_1u}_{C^1}$.
So
$\norm{d_2d_1u}^2_{C^2}\leq\frac{4|T|^2\e^2}{\delta^4}\norm{d_1u}_{C^1}$.
Since $d_1(B^0(\frac{\delta^2}{|T|}))$ is dense in
$B^1(\frac{\delta^2}{|T|})$, for any $f\in
B^1(\frac{\delta^2}{|T|})$ it also holds that
$\norm{d_2f}^2_{C^2}\leq\frac{4|T|^2\e^2}{\delta^4}\norm{f}_{C^1}$.
 \end{proof}

 \begin{cor}\label{C1}
For every $f\in C^1$ one has
 $$
\langle
d^*_1f,d^*_1f\rangle_{C^0}\geq\Bigl(4-\frac{2}{\lambda_1}-\frac{20\e}{3\lambda_1}-\frac{8|T|^2\e^2}{3\lambda_1\delta^4}\Bigr)\langle
f,f\rangle_{C^1}.
 $$
 \end{cor}

Corollary \ref{C1} provides the constant $c$ for Proposition
\ref{tm1}. Now the function $\alpha(\varepsilon)$ from Theorem
\ref{1} should satisfy
$\alpha(\varepsilon)=\max\bigl\lbrace\frac{10\varepsilon}{3\lambda_1}+\frac{8|T|^2\varepsilon^2}{3\lambda_1\delta^4},\delta\bigr\rbrace$.
In order to get a continuous function with $\alpha(0)=0$ one may
take $\delta=\varepsilon^{2/5}$. Then
$\alpha(\varepsilon)=O(\varepsilon^{2/5})$ and Theorem \ref{1}
directly follows from Corollary \ref{C1} and Proposition
\ref{tm1}. \qed

\bigskip
The authors are grateful to N. Monod for useful remarks.

\end{document}